\documentclass[a4paper,12pt]{amsart}
\usepackage{amssymb}
\usepackage{ifthen}
\usepackage{graphicx}

\newtheorem{thm}{Theorem}
\newtheorem{cor}{Corollary}
\newtheorem{lem}{Lemma}

\newtheorem{rem}{Remark}

\newtheorem{conj}{Conjecture}
\theoremstyle{definition}
\newtheorem{example}[equation]{Example}
\newtheorem{prob}[equation]{Problem}

\newcounter {own}
\def\theown {\thesection       .\arabic{own}}

\newenvironment{pf}[1][]{%
 \vskip 3mm
 \noindent
 \ifthenelse{\equal{#1}{}}%
  {{\slshape Proof. }}%
  {{\slshape #1.} }%
 }%
{\qed\bigskip}

\newcounter{alphabet}
\newcounter{tmp}
\newenvironment{Thm}[1][]{\refstepcounter{alphabet}%
\bigskip%
\noindent%
{\bf Theorem \Alph{alphabet}}%
\ifthenelse{\equal{#1}{}}{}{ (#1)}%
{\bf .} \itshape}{\vskip 8pt}

\makeatletter
\newcommand{\Ref}[1]{\@ifundefined{r@#1}{}{\setcounter{tmp}{\ref{#1}}\Alph{tmp}}}
\makeatother

\newcommand{\IN}{{\mathbb N}}
\newcommand{\IC}{{\mathbb C}}
\newcommand{\ID}{{\mathbb D}}

\def\be{\begin{equation}}
\def\ee{\end{equation}}

\newcommand{\bee}{\begin{enumerate}}
\newcommand{\eee}{\end{enumerate}}

\newcommand{\blem}{\begin{lem}}
\newcommand{\elem}{\end{lem}}
\newcommand{\bthm}{\begin{thm}}
\newcommand{\ethm}{\end{thm}}
\newcommand{\bcor}{\begin{cor}}
\newcommand{\ecor}{\end{cor}}
\newcommand{\beg}{\begin{example}}
\newcommand{\eeg}{\end{example}}
\newcommand{\begs}{\begin{examples}}
\newcommand{\eegs}{\end{examples}}
\newcommand{\bdefe}{\begin{defin}}
\newcommand{\edefe}{\end{defin}}
\newcommand{\bprob}{\begin{prob}}
\newcommand{\eprob}{\end{prob}}
\newcommand{\bei}{\begin{itemize}}
\newcommand{\eei}{\end{itemize}}

\newcommand{\bcon}{\begin{conj}}
\newcommand{\econ}{\end{conj}}
\newcommand{\bcons}{\begin{conjs}}
\newcommand{\econs}{\end{conjs}}
\newcommand{\bprop}{\begin{propo}}
\newcommand{\eprop}{\end{propo}}
\newcommand{\br}{\begin{rem}}
\newcommand{\er}{\end{rem}}
\newcommand{\brs}{\begin{rems}}
\newcommand{\ers}{\end{rems}}
\newcommand{\bo}{\begin{obser}}
\newcommand{\eo}{\end{obser}}
\newcommand{\bos}{\begin{obsers}}
\newcommand{\eos}{\end{obsers}}
\newcommand{\bpf}{\begin{pf}}
\newcommand{\epf}{\end{pf}}
\newcommand{\ba}{\begin{array}}
\newcommand{\ea}{\end{array}}
\newcommand{\beq}{\begin{eqnarray}}
\newcommand{\beqq}{\begin{eqnarray*}}
\newcommand{\eeq}{\end{eqnarray}}
\newcommand{\eeqq}{\end{eqnarray*}}

\newcommand{\Llra}{\Longleftrightarrow}

\newcommand{\ra}{\rightarrow}

\newcommand{\ds}{\displaystyle}

\def\cc{\setcounter{equation}{0}   
\setcounter{figure}{0}\setcounter{table}{0}}

\newcounter{minutes}
\divide\time by 60
\newcounter{hours}
\multiply\time by 60 \addtocounter{minutes}{-\time}

\begin{document}
\bibliographystyle{amsplain}
\title[Univalence of the average of two analytic functions]{Univalence of the average of two analytic functions}

\thanks{
File:~\jobname .tex,
          printed: \number\day-\number\month-\number\year,
          \thehours.\ifnum\theminutes<10{0}\fi\theminutes}


\author{M. Obradovi\'{c}}
\address{M. Obradovi\'{c}, Department of Mathematics,
Faculty of Civil Engineering, University of Belgrade,
Bulevar Kralja Aleksandra 73, 11000
Belgrade, Serbia. } \email{obrad@grf.bg.ac.rs}

\author{S. Ponnusamy ${}^\dagger$}
\address{S. Ponnusamy, Department of Mathematics,
Indian Institute of Technology Madras, Chennai--600 036, India.}
\email{samy@iitm.ac.in}

\keywords{Coefficient inequality, partial sums, radius of univalence, analytic,
univalent and starlike functions\\
$^\dagger$ Corresponding author}

\begin{abstract}
Let $\mathcal{A}$ denote the set of all analytic functions $f$ in the unit disk
$\ID=\{z:\,|z|<1\}$ of the form $f(z)=z+\sum_{n=2}^{\infty}a_nz^n.$
Let $\mathcal{U}$ denote the set of all $f\in \mathcal{A}$, $f(z)/z\neq 0$ and
satisfying the condition
$$\left | f'(z)\left (\frac{z}{f(z)} \right )^{2}-1\right | < 1
~\mbox{ for $z\in \ID$}.
$$
Functions in ${\mathcal U}$ are known to be univalent in $\ID$. For $\alpha \in [0,1]$, let
$$ \mathcal{N}(\alpha )=\left \{f_\alpha :\,
f_\alpha (z)=(1-\alpha )f(z)+\alpha \int_0^z\frac{f(t)}{t}\,dt,~ \mbox{ $f\in\mathcal{A}$ with $|a_n|\leq n$ for
$n\geq 2$}\right\}.
$$
In this paper, we first show that the condition $\sum_{n=2}^{\infty}n|a_n|\leq 1$ is sufficient for $f$
to be in ${\mathcal U}$ and the same condition is necessary for $f\in {\mathcal U}$ in case all $a_n$'s are
negative. Next, we obtain the radius of  univalence of functions in the class $\mathcal{N}(\alpha )$.
Also, for  $f,g\in \mathcal{U}$ with $\frac{f(z)+g(z)}{z}\neq 0$ in $\ID$,
$F(z)=(f(z)+g(z))/2$, and $G(z)=r^{-1}F(rz)$, we determine a range of $r$ such that
$G\in {\mathcal U}$. As a consequence of these results, several special cases are presented.
\end{abstract}
\thanks{The work of the first author was supported by MNZZS Grant, No. ON174017, Serbia.}

\maketitle
\pagestyle{myheadings}
\markboth{M. Obradovi\'{c}, and S.Ponnusamy}{Average of two univalent functions}
\cc

\section{Introduction}
Let $\mathcal A$ denote the set of all
normalized analytic functions $f$ in the unit disk $\ID=\{z\in\IC:\, |z|<1\}$ of the form
$f(z)=z+\sum_{n=2}^{\infty}a_nz^n.$
For $f\in\mathcal A$ and $n\in \IN$, we write
$s_{n}(f) (z)=z+\sum_{k=2}^{n}a_kz^k,
$
for the $n$-th partial sums or sections of $f$.
Also, we denote by ${\mathcal S}$ the class of all univalent functions in $\mathcal A$. It is
well-known that if $f\in {\mathcal S} $ then $|a_n|\leq n$ for $n\geq 2$.
The class of convex and starlike functions are two important geometric subclasses of
$ \mathcal{S}$, denoted by $\mathcal C$ and ${\mathcal S}^*$, respectively (see \cite{Du,Go}).
For $f \in {\mathcal C}$, one has $|a_n|\leq 1$ for all $n\ge 2$.
Two other subclasses of $\mathcal S$ that are studied extensively are
\beqq
{\mathcal R}_1 & =& \left \{f\in {\mathcal A}:\, |f'(z)-1|<1  ~\mbox{ for $z\in \ID$}\right \}\\
{\mathcal S}_1^* & =& \left \{f\in {\mathcal S}:\, |(zf'(z)/f(z))-1|<1  ~\mbox{ for $z\in \ID$}\right \}.
\eeqq
Finally, we define (\cite{Ak,OP-01,obpo-2007a})
$$\mathcal{U}=\left \{f\in {\mathcal A}:\, \mbox{$f(z)/z\neq 0$ and }~ \left | U_f(z) -1\right | < 1 ~\mbox{ for
$z\in \ID$}\right \}, 
$$
where
$$U_f(z)= f'(z)\left (\frac{z}{f(z)} \right )^{2}.
$$
In the recent years, the class $\mathcal{U}$ and its association with a number of subclasses of $\mathcal{S}$
together with certain integral transformations have been studied in details (see \cite{FR-2006,OP-01,obpo-2007a,OPSV}).
It is well-known that $\mathcal{U}\subsetneq\mathcal{S}$ (see \cite{Ak,OzNu}). It is interesting to observe that
the Koebe function $k(z)=z/(1-z)^2$ belongs to $\mathcal{U}$.
The following result is well-known.

\begin{Thm}
{\rm \cite[Theorem 3]{Cl-Keo1960}} \label{6-10thA}
If $\sum_{n=2}^{\infty}n|a_n|\leq 1,$ then $f(z)=z+\sum_{n=2}^{\infty}a_nz^n$
belongs to ${\mathcal S}_1^*\cap{\mathcal R}_1$.
\end{Thm}

Let $\mathcal F$ and $\mathcal G$ be two subclasses of $\mathcal A$.
If for every $f\in {\mathcal F}$, $G_r(z)=r^{-1} f(rz)\in\mathcal{G}$ for $r\leq r_0$, and
$r_0$ is the maximum value for which this holds, then we say that $r_0$ is
the $\mathcal{G}$-radius in $\mathcal{F}$. That is,
$$\ds r_0=r_{\mathcal{G}}(\mathcal{F}) :=\sup\left\{r>0: \,G_r\in\mathcal{G} \right\}
~\mbox{ for every $f\in\mathcal{F}.$}
$$
There are many results of this type in
the theory of univalent functions, see \cite[Chapter~13, Vol.~2]{Go}. For instance, the ${\mathcal C}$-radius
in ${\mathcal S}$ is $2-\sqrt{3}$ which is referred to as the \textit{radius of convexity} of ${\mathcal S}$.
Similarly, the \textit{radius of starlikeness} for the class $\mathcal S$ is $\tanh (\pi/4) \approx 0.65579$
(see \cite[p.~44]{Du} and \cite[Chapter~8, Vol.~1]{Go}).  We now recall the following result due to
Gavrilov \cite[Theorem 1]{Gav-1970} which does not seem to be known for many readers.

\begin{Thm}\label{6-10-ThB}
Let $\mathcal{F}=\{f\in \mathcal{A}:\, |a_{n}|\leq n \mbox{ for $n\geq 2$}\}.$
Then $f$ and each of its $n$-th partial sum $s_n(f)$ is univalent for $|z|<r_{\mathcal{S}}$, where
$r_{\mathcal{S}}\approx 0.164878$ is the root of the  equation $2(1-r)^{3}-r-1=0$ in $(0,1)$. Here $r_{\mathcal{S}}$
is the radius of univalence of
$$g(z)=2z -\frac{z}{(1-z)^{2}}, ~~z\in \ID.
$$
Equivalently,
$r_{\mathcal{S}}$ is the $\mathcal{S}$-radius in $\mathcal{F}$.
\end{Thm}

In the forgoing discussion, we say that  $f$ belongs to ${\mathcal U}$ in the disk $|z|<r$
if the inequality in the above definition of ${\mathcal U}$ holds for $|z|< r$ instead of the whole unit disk $\ID$.
In other words, this is equivalent to saying that $G$ defined by $G(z)=r^{-1}f(rz)$ belongs to $\mathcal{U}$,
when $f$ belongs to ${\mathcal U}$ in the disk $|z|<r$.

In order to indicate another example of the type as in Theorem \Ref{6-10-ThB}, we may recall the following result
which is indeed a corollary to a general result (see \cite{OP-2005}).

\begin{Thm}\label{cor-radi1}
$r_{\mathcal{U}}(\mathcal{S})= 1/\sqrt{2}$.
\end{Thm}

The paper is organized as follows. The main results are stated in Section \ref{6-10sec2}
and their proofs and some of their consequences are presented in Section \ref{6-10sec3}.
First we present a sufficient coefficient condition for a
function $f$ to be in $\mathcal{U}$. The condition is also shown to be necessary if the
coefficients of $f\in {\mathcal A}$ are negative. Next, we present several simple observations
concerning the radius of univalence (and of starlikeness) of certain class of analytic functions
in the unit disk. Finally, we obtain a radius of univalence of average of two
functions from $\mathcal{U}$.

\section{Main Results}\label{6-10sec2}

We now state our first result.

\bthm \label{6-10th4}
If $\sum_{n=2}^{\infty}n|a_n|\leq 1,$ then $f(z)=z+\sum_{n=2}^{\infty}a_nz^n$
belongs to ${\mathcal U}$. The result is sharp.
\ethm

At this place, it is worth recalling that the class $\mathcal{U}$ is
neither included in $\mathcal{S}^*$ nor includes the class ${\mathcal S}^*$.  Also, the class $\mathcal{U}$ is neither
contained in $\mathcal{R}_1$ nor contains the class ${\mathcal R}_1$ (see for example \cite{obpo-2007a}).
For instance, the function $f$ defined by
$$f(z)=\frac{z}{1+\frac{1}{2}z+\frac{1}{2}z^3}
$$
is in ${\mathcal U}\backslash {\mathcal S}^*$ (see also \cite{FR-2006}). Indeed for this function
$$\frac{zf'(z)}{f(z)}=\frac{1-z^3}{1+\frac{1}{2}z+\frac{1}{2}z^3}
$$
we have at $z_0= (-1+i)/\sqrt 2$,
$$\frac{z_0f'(z_0)}{f(z_0)} =\frac{2-2\sqrt{2}}{3} +\frac{1-2\sqrt 2}{3}i
$$
which gives that ${\rm Re}\, \{z_0f'(z_0)/f(z_0)\}<0$ showing that $f$ is not starlike in $\ID$.
Moreover, although the Koebe function $k(z)$ belongs
to ${\mathcal U}\cap {\mathcal S}^*$, it is known that ${\mathcal S}\not\subset  {\mathcal U}$.
In view of these observations, Theorem \ref{6-10th4} refines Theorem \Ref{6-10thA} and hence,
Theorem \ref{6-10th4} will be useful to generate nice class of examples of functions in ${\mathcal U}$
that are also starlike in $\ID$.

For functions with negative coefficients, the converse of Theorem \ref{6-10th4} is also true.

\bthm \label{6-10th5}
A function  $f(z)=z-\sum_{n=2}^{\infty}|a_n|z^n$ is in the class ${\mathcal U}$ if and only if
$\sum_{n=2}^{\infty}n|a_n|\leq 1$.
\ethm

From the result of Silverman \cite{Silv-1975}, we may now formulate the above discussion as

\bcor
Suppose that $f(z)=z-\sum_{n=2}^{\infty}|a_n|z^n$ belongs to $\mathcal A$. Then we have the following equivalent statements:
$$f\in {\mathcal U} \Llra f\in {\mathcal S}^*_1 \Llra f\in {\mathcal R}_1
\Llra \sum_{n=2}^{\infty}n|a_n|\leq 1.
$$
\ecor

In connection with a problem due to \cite{Hayman-67}, Ruscheweyh and Wirths \cite{Rus-Wirth76}
discussed the univalency of functions in the set of convex linear combinations of the form
$$ \mu f(z) +(1-\mu)g(z), ~\mu \in [0,1],
$$
when $f,g$ belonging to suitable subsets of $\mathcal S$. We shall consider a similar problem
which is indeed a generalization of Theorem \Ref{6-10-ThB}. For $\alpha \in [0,1]$, let
$$ \mathcal{N}(\alpha )=\left \{f_\alpha :\,
f_\alpha (z)=(1-\alpha )f(z)+\alpha \int_0^z\frac{f(t)}{t}\,dt,~ \mbox{ $f\in\mathcal{A}$ with $|a_n|\leq n$ for
$n\geq 2$}\right\}
$$
and for convenience, we set $r_{\mathcal{S}}(\alpha):=r_{\mathcal{S}}({\mathcal N}(\alpha))$.

\bthm\label{6-10th1}
The number $r_{\mathcal{S}}(\alpha)$ is the root in $(0,1)$ of the  equation
\be\label{6-10eq3}
2(1-r)^{3} +(2\alpha -1)r-1=0.
\ee
The extremal function is
$$g(z)=2z-\frac{\alpha z}{1-z} -\frac{(1-\alpha)z}{(1-z)^{2}}, ~~z\in \ID.
$$
\ethm

For a ready reference, the values of $r_{\mathcal{S}}(\alpha)$
for certain values of $\alpha \in [0,1]$ are listed in Table \ref{table1}.
\begin{table}
\center{
\begin{tabular}{|c|c|}
\hline
$\alpha $ & $r_{\mathcal S}(\alpha)$
\\
\hline
0 & 0.164878
\\
$1/5$& 0.178866
\\
$1/4$ & 0.182817
\\
$1/3$ & 0.189894
\\
$1/2$ & 0.206299
\\
$3/4$ & 0.23931
\\
$4/5$& 0.247756
\\
1 & 0.292893

\\
\hline
\end{tabular}
}
\bigskip
\caption
{The radius $r_{\mathcal{S}}(\alpha)$ of univalence for $f_\alpha\in {\mathcal N}(\alpha)$\label{table1}}
\end{table}

\br
{\rm
If $f\in {\mathcal S}$, then the Biernacki integral $(Bf)(z)=\int_0^z(f(t)/t)\,dt$ is not necessarily
univalent in $\ID$.
}
\er

The case $\alpha =0$ leads to Theorem \Ref{6-10-ThB}.
Set $\mathcal{N}=\mathcal{N}(1)$ so that
$$\mathcal{N}=\{f\in \mathcal{A}:\, |a_{n}|\leq 1 \mbox{ for $n\geq 2$}\}
$$
and therefore, the case $\alpha =1$ of Theorem \ref{6-10th1} may be reformulated as

\bcor\label{6-10cor2}
Then $f\in \mathcal{N}$ and each of its $n$-th partial sum $s_n(f)$ is univalent for $|z|<r$, where
$r\approx 1-\frac{\sqrt{2}}{2}\approx 0.292893$ is the  root of the  equation
$2(1-r)^{3}+r-1=0$ in $(0,1)$. The extremal function is
$$g(z)=2z-\frac{z}{1-z}, ~~z\in \ID.
$$
\ecor
We remark that
$$ \mathcal{N} \supset \{f\in{\mathcal A}:\, {\rm Re}\, (f(z)/z)>1/2,~z\in \ID\} \supset \{f:\, f\in {\mathcal C}\}.
$$

In the case of $\alpha =1/2$, we see that if
$$ \mathcal{N}(1/2)=\left \{f_{1/2}(z)=\frac{1}{2} \left (f(z)+ \int_0^z\frac{f(t)}{t}\,dt\right ):\, f\in {\mathcal A}~
\mbox{ and  $|a_n|\leq n$ for $n\geq 2$}\right\},
$$
then functions $f_{1/2}$ in $\mathcal{N}(1/2)$ are univalent in $|z|<1-2^{-1/3}\approx 0.206299$, but not necessarily
in any larger disk as the function
$$g(z)=2z-\frac{z}{2(1-z)} -\frac{z}{2(1-z)^{2}}, ~~z\in \ID,
$$
shows.

In order to motivate our next theorem, let us consider the sum of two univalent functions $f,g\in \mathcal{S}$.
Then the average function $F$ defined by
$$F(z)=\frac{f(z)+g(z)}{2}, ~~z\in \ID,
$$
belongs to the class $\mathcal{N}(0)$, and hence, by Theorem \Ref{6-10-ThB},
we conclude that the radius of univalence of the function $F$ is not
smaller than that of the number given in Theorem \Ref{6-10-ThB}.

Similarly, if $f,g\in \mathcal{C}$ then the necessary coefficient estimates on the Taylor coefficients of $f,g$
show that the average function $F$
defined by
$$F(z)=\frac{f(z)+g(z)}{2}, ~~z\in \ID,
$$
belongs to the class $\mathcal{N}$ defined in Corollary \ref{6-10cor2}. Again,
the radius of univalence of the function $F$ is not
smaller than that given in Corollary \ref{6-10cor2}. Thus, any result which
gives radius bigger than that given in these two corollaries will provide us results
with an improved bound for the radius. Our next result fills this idea
if we consider the subclass $\mathcal U$ of $\mathcal S$ for our investigation.

In the following theorem, we use certain well-known  basic facts.
Each $f\in \mathcal{U}$ with $a_2=f''(0)/2$ can be written in the form
\be\label{eq1}
f'(z)\left (\frac{z}{f(z)}\right )^{2}= -z\left(\frac{z}{f(z)}\right)'+\frac{z}{f(z)}= 1+w(z), ~~z\in \ID,
\ee
where $w\colon \ID \ra \ID$ is analytic with $w(0)=w'(0)=0$. Consequently, the classical Schwarz lemma gives
$$\left|U_f(z)-1\right|\leq |z|^{2}~\mbox{ for $z\in \ID$}
$$
and, by (\ref{eq1}), we easily have
\be\label{8veq1a}
\frac{z}{f(z)}=1 -a_2z-\int_0^1\frac{w(tz)}{t^2}\,dt,  \quad z\in \ID
\ee
so that
\be\label{6-10eq9}
\left |\frac{z}{f(z)}-1 -a_2z\right |\leq |z|^2,  \quad z\in \ID.
\ee

\bthm\label{6-10th2}
Let $f,g\in \mathcal{U}$ with $\frac{f(z)+g(z)}{z}\neq 0$ in the unit disk $\ID$, and
$F$ be defined by
\be\label{6-10eq4}
F(z)=\frac{f(z)+g(z)}{2}, ~~z\in \ID.
\ee
Then $F\in {\mathcal U}$ for $|z|<r_0$, where $ r_{0}\approx 0.262453 $ is the root of the equation
\be\label{6-10eq8}
r^{2}+(1+r^{2})\frac{(2r+r^{2})^{2}}{\left(1-2r-\frac{\pi}{\sqrt[4]{90}}\frac{r^{2}}{\sqrt[4]{1-r^{4}}}\right)^{2}}=1.
\ee
\ethm

Using Theorem \Ref{cor-radi1} and Theorem \ref{6-10th2}, one can quickly deduce the following result.

\bthm\label{6-10th3}
Let $f,g\in \mathcal{S}$ with $\frac{f(z)+g(z)}{z}\neq 0$ in the unit disk $\ID$.
Then the function $F$ define by {\rm (\ref{6-10eq4})} belongs to ${\mathcal U}$
for $|z|<\frac{1}{\sqrt{2}}r_{0}$, where $ r_{0}\approx 0.262453 $ is the root of the equation
{\rm (\ref{6-10eq8})}.
\ethm

At this place it worth recalling that $\mathcal{S}\subset \mathcal{N}(0)$.

\section{A Lemma and an Example}\label{6-10sec3}
For the proof of Theorem \ref{6-10th2} and the discussion in Example~1 below, 
we need the following lemma.

\blem\label{5-10lem1}
Let $\phi(z)=1+\sum_{n=1}^\infty b_nz^n$ be a non-vanishing analytic
function on $\ID$ and let $f$ be of the form 
$$f(z)=\frac{z}{\phi (z)}.
$$
Then, we have the following:
\begin{enumerate}
\item[\textbf{(a)}] If
$\sum_{n=2}^\infty (n-1)|b_n|\leq 1 ,$ then $f\in {\mathcal U} $.
\item[\textbf{(b)}] If $\sum_{n=2}^\infty (n-1)|b_n|\leq 1-|b_1|$, then $f\in {\mathcal S}^*$.
\item[\textbf{(c)}] If $f\in {\mathcal U}$, then $\sum_{n=2}^{\infty}(n-1)^2|b_n|^2\leq 1.$
\end{enumerate}
\elem

The conclusion \textbf{(a)} in Lemma \ref{5-10lem1} is from \cite{OP-01,OPSV} whereas the
\textbf{(b)} is due to Reade et. el. \cite[Theorem 1]{RST-84}. Finally, as  $f\in {\mathcal U}$, we have
$$\left |f'(z)\left (\frac{z}{f(z)} \right )^{2}-1\right |=
\left |-z \left ( \frac{z}{f(z)} \right )'+\frac{z}{f(z)} -1\right | = \left |\sum_{n=2}^{\infty}(n-1)b_{n}z^{n}\right |
\leq 1
$$
and so \textbf{(c)} is an immediate consequence of Gronwall's area theorem, and can be obtained also from Parseval's
relation.

\vspace{8pt}

\noindent
{\bf Example 1.~} Consider
$$f(z)=\frac{z}{(1-iz)^{2}} ~\mbox{ and }~g(z)=\frac{z}{(1+iz)^{2}}, ~~z\in \ID.
$$
Then $f,g\in\mathcal S$, and for $\alpha \in [0,1]$, we obtain
$$f_\alpha(z)=\frac{(1-\alpha)z}{(1-iz)^{2}}+\frac{\alpha z}{1-iz}
 ~\mbox{ and }~g_\alpha(z)=\frac{(1-\alpha)z}{(1+iz)^{2}}+\frac{\alpha z}{1+iz}, ~~z\in \ID.
$$
We observe that  $f_\alpha, g_\alpha\in \mathcal{N}(\alpha )$. Now, we introduce
$$F_\alpha(z)=\frac{1}{2}(f_\alpha(z)+g_\alpha(z)), ~~z\in \ID.
$$
Simple calculation shows that
$$F_\alpha(z) 
=\frac{z[1-(1-2\alpha)z^2]}{(1+z^2)^2}
$$
and
$$\frac{z}{F_\alpha (z)}=\frac{(1+z^{2})^{2}}{1-(1-2\alpha)z^2}=1+(3-2\alpha)z^2
+\frac{4(1-\alpha)^2z^4}{1-(1-2\alpha)z^2}, ~~z\in \ID.
$$
For $0<r\leq1$, define $G_\alpha$ by $G_{\alpha, r}(z)=r^{-1}F_\alpha(rz)$ so that
$$\frac{z}{G_{\alpha, r}(z)}=1+\sum_{n=2}^{\infty}B_{2n}(r)z^{2n}
$$
where
$$B_{2n}(r) =\left\{\ba{rl}
(3-2\alpha)r^2 &\mbox{ if $n=1$}\\
4(1-\alpha)^2(1-2\alpha)^{n-2}r^{2n} &\mbox{ if $n\ge 2$}.
\ea\right.
$$
By Lemmas \ref{5-10lem1}\textbf{(a)} and \textbf{(b)}, the function $G_{\alpha, r} \in
{\mathcal U}\cap {\mathcal S}^*$ if $r$ satisfies the condiiton
$$S(r)=\sum_{n=2}^\infty (2n-1)|B_{2n}(r)|\leq 1.
$$
A computation gives
\beqq
S(r) &=& (3-2\alpha)r^2 +4(1-\alpha)^2r^4\sum_{n=2}^{\infty}(2n-1)|1-2\alpha|^{n-2} r^{2(n-2)}\\
&=& (3-2\alpha)r^2 +4(1-\alpha)^2r^4\left \{\frac{3-|1-2\alpha|r^2}{(1-|1-2\alpha|r^2)^2}
\right \}.
\eeqq
Thus, $S(r)\leq 1$ if and only if $r$ satisfies the inequality
$$ (3-2\alpha)r^2 (1-|1-2\alpha|r^2)^2 + 4(1-\alpha)^2r^4(3-|1-2\alpha|r^2) -(1-|1-2\alpha|r^2)^2\leq 0.
$$
If $1-2\alpha >0$, then the last condition is equivalent to
$$(1+r^2) \left ((1-2\alpha)r^4-6(1-\alpha)r^2+1\right )\geq 0
$$
which is true if $r^2 \leq K_\alpha ^2$,
where
$$K_\alpha =\sqrt{\frac{3(1-\alpha)- \sqrt{9\alpha ^2-16\alpha +8}}{1-2\alpha}} ~\mbox{ for $0\leq \alpha <1/2$}.
$$
Thus, for $0\leq \alpha  <1/2$, the inequality $S(r)\leq 1$ holds if $r\leq K_\alpha $. Moreover, for $\alpha =1/2$,
we obtain that $S(r)\leq 1$ whenever
$$r\leq \lim_{\alpha \ra (1/2)-}K_\alpha =\frac{1}{\sqrt{3}}.
$$

Finally, if $1-2\alpha <0$ then the inequality $S(r)\leq 1$ holds if and only if
$$(3-2\alpha)r^2 (1+(1-2\alpha)r^2)^2 + 4(1-\alpha)^2r^4(3+ (1-2\alpha)r^2) -(1+(1-2\alpha)r^2)^2\leq 0,
$$
or equivalently $r\leq K_\alpha'$, where $K_\alpha'$ for $1/2<\alpha \leq 1$ is
the root of the equation
$$-1+ (1+2\alpha)r^2 +(17-36\alpha +16\alpha ^2)r^4 + (7-16\alpha +8\alpha^2)(1-2\alpha)r^6=0
$$
in the unit interval $(0,1)$. Setting
$$ r(\alpha)=\left \{ \ba{cl}
K_\alpha & \mbox{ if $0\leq \alpha  <1/2$} \\
\ds \frac{1}{\sqrt{3}} &  \mbox{ if $\alpha =1/2$} \\
K_\alpha'  &\mbox{ if $1/2<\alpha \leq 1$},
\ea
\right .
$$
we see that
\begin{table}
\center{
\begin{tabular}{|c|c|}
\hline
$\alpha $ & $r(\alpha) 
$
\\
\hline
0 & 0.414214
\\
$1/5$& 0.462667
\\
$1/4$ & 0.477491
\\
$1/3$ & 0.505408
\\
$1/2$ & 0.57735
\\
$3/4$ & 0.69398
\\
$4/5$& 0.727725
\\
1 & 1
\\
\hline
\end{tabular}
}
\bigskip
\caption
{Radius of univalence and starlikeness of $F_\alpha(z)$. \label{table2}
}
\end{table}
the function $G_{\alpha, r}(z)$ is univalent in $\ID$ if $0<r\leq r(\alpha)$. Equivalently, it means
that the function $F_\alpha$ is univalent in the disk $|z|<r(\alpha)$ and the result is sharp.
In Table \ref{table2}, we list the values of $r(\alpha)$ for certain choices of $\alpha$.

The above discussion suggests the following

\vspace{8pt}
\noindent
{\bf Problem 1.} Suppose that for $\alpha \in [0,1]$,
$$ \mathcal{S}(\alpha )=\left \{
f_\alpha (z)=(1-\alpha )f(z)+\alpha \int_0^z\frac{f(t)}{t}\,dt:\, f\in\mathcal{S}\right\}.
$$
and
$$ \mathcal{F}(\alpha )=\left\{F_\alpha(z)=\frac{1}{2}(f_\alpha(z)+g_\alpha(z)):\, f_\alpha ,g_\alpha \in\mathcal{S}(\alpha )\right \}.
$$
Determine the radius of univalence of functions in $F_\alpha\in\mathcal{F}(\alpha )$.

\section{Proofs and some of their consequences 
}\label{6-10sec5}
\noindent \textbf{Proof of Theorem \ref{6-10th4}.} Suppose that $\sum_{n=2}^{\infty}n|a_n|\leq 1$
and $f(z)=z+\sum_{n=2}^{\infty}a_nz^n$.
In order to show that $f\in{\mathcal U}$, we need to show that $| U_f(z) -1| < 1$ for
$z\in \ID$.
As
\beqq
\left | f'(z)-\left (\frac{f(z)}{z} \right )^2\right |
&=& \left | 1 + \sum_{n=2}^{\infty}na_nz^{n-1}-\left (1+\sum_{n=2}^{\infty}a_nz^{n-1} \right )^2\right | \\
&=& \left | \sum_{n=2}^{\infty}(n-2)a_nz^{n-1}-\left (\sum_{n=2}^{\infty}a_nz^{n-1} \right )^2\right | \\
&=& |z|^2\left | \sum_{n=3}^{\infty}(n-2)a_nz^{n-3}-\left (\sum_{n=2}^{\infty}a_nz^{n-2} \right )^2\right |,
\eeqq
it follows that
\beqq
\left | f'(z)-\left (\frac{f(z)}{z} \right )^2\right |
& < & \sum_{n=2}^{\infty}(n-2)|a_n| +\left (\sum_{n=2}^{\infty}|a_n|\right )^2 \\
& \leq & 1-2\sum_{n=2}^{\infty}|a_n| +\left (\sum_{n=2}^{\infty}|a_n|\right )^2 \\
& \leq  & \left (1-\sum_{n=2}^{\infty}|a_n|\right )^2
\leq  \left | \frac{f(z)}{z} \right |^2
\eeqq
which implies that
$$ | U_f(z) -1| =\left |f'(z)\left (\frac{z}{f(z)} \right )^{2}-1\right |< 1, \quad z\in \ID.
$$
and hence, $f\in{\mathcal U}$.

To see that the upper bound $1$ in the coefficient condition cannot be replaced by $1+\epsilon$, $\epsilon >0$, we
consider the function
$$f_\epsilon (z)=z+\frac{1+\epsilon}{n}z^n  ~~(n\geq 2).
$$
We obtain that
$f_\epsilon '(z)=1+(1+\epsilon)z^{n-1}$ has a zero in $\ID$ as $\epsilon >0$ and hence, $f_\epsilon$ is not univalent in $\ID$.
In particular, $f_\epsilon\notin {\mathcal U}$. Thus, the result is sharp.
\hfill $\Box$

\vspace{8pt}

\br
{\rm
Theorem \Ref{6-10thA} in particular gives condition on the Taylor coefficients of $f$ so that
$f$ and its $n$-th partial sum $s_n(f)$ is not only univalent but is starlike in $\ID$.
In view of this observation, from the proof of
Theorem \ref{6-10th1}, it follows that the quantity $r_{\mathcal{S}}(\alpha)$ in Theorem \ref{6-10th1}
is indeed the $\mathcal{S}$-radius in $\mathcal{N}(\alpha )$ as well as
the radius of starlikeness of functions in the class $\mathcal{N}(\alpha)$. These observations
and the proof of Theorem \ref{6-10th4} show that Theorem \ref{6-10th4} may
be stated in an improved form.
}
\er

\noindent \textbf{Proof of Theorem \ref{6-10th5}.}
In view of Theorem \ref{6-10th4}, it suffices to show the only if part. Assume
that $|U_f(z) -1|<1$ in $\ID$. That is,
\beqq
\left | f'(z)\left (\frac{z}{f(z)} \right )^{2} -1\right | &=&  \left | \frac{\ds 1 - \sum_{n=2}^{\infty}n|a_n|z^{n-1}}
{\ds \left (1-\sum_{n=2}^{\infty}|a_n|z^{n-1} \right )^2} -1\right | \\
&=& |z|^2\left | \frac{\ds - \sum_{n=3}^{\infty}(n-2)|a_n|z^{n-3} -\left (\sum_{n=2}^{\infty}|a_n|z^{n-2} \right )^2}
{ \left (\ds 1-\sum_{n=2}^{\infty}|a_n|z^{n-1} \right )^2} \right | <1
\eeqq
for $z\in \ID$. Choose values of $z$ on the real axis so that $U_f(z) -1$ is real. Upon clearing the
denominator in the last expression and letting $z\rightarrow 1^{-}$ through real values, we obtain
$$\sum_{n=3}^{\infty}(n-2)|a_n| +\left (\sum_{n=2}^{\infty}|a_n|\right )^2
\leq  \left (\ds 1-\sum_{n=2}^{\infty}|a_n| \right )^2.
$$
Thus, $\sum_{n=2}^{\infty}n|a_n|\leq 1$, and the proof is complete.
\hfill $\Box$

\vspace{8pt}

\noindent \textbf{Proof of Theorem \ref{6-10th1}.}
Let $f_\alpha \in \mathcal{N}(\alpha )$ with $f(z)=z+\sum_{n=2}^{\infty}a_nz^n$. Then $f_\alpha'$ takes the form
$$f_\alpha '(z)= (1-\alpha )f'(z)+\alpha \frac{f(z)}{z}=
1+\sum_{n=2}^{\infty}((1-\alpha )n+\alpha) a_nz^{n-1}.
$$
As $|a_n|\leq n$ for $n\geq 2$, it follows that for $|z|=r$
\beqq
|f_\alpha '(z)-1|&\leq & \sum_{n=2}^{\infty}((1-\alpha )n+\alpha )nr^{n-1}\\
&=& (1-\alpha )\left (\frac{1+r}{(1-r)^{3}}-1\right )
+\alpha\left (\frac{1}{(1-r)^{2}}-1\right )
\eeqq
so that $|f_\alpha '(z)-1|<1$ whenever $R_\alpha (r)>0$, where
$$R_\alpha (r)=2(1-r)^{3} +(2\alpha -1)r-1 .
$$
Consequently, by (\ref{6-10eq3}), ${\rm Re}\, f_\alpha '(z)>0 $ for $|z|<r_0$ where
$r_0\geq r_{\mathcal{S}}(\alpha)$. Next we show that $r_0=r_{\mathcal{S}}(\alpha)$. For this,
we consider the function
$$g(z)=z-\sum_{n=2}^{\infty}((1-\alpha )n+\alpha)z^{n}\in \mathcal{N}(\alpha ).
$$
It is a simple exercise to see that
$$g(z)=\frac{z(1-(4-\alpha)z+2z^{2})}{(1-z)^{2}} =2z-\frac{\alpha z}{1-z} -\frac{(1-\alpha)z}{(1-z)^{2}}
$$
so that
$$g'(z)=\frac{2(1-z)^3-1 +(2\alpha-1)z}{(1-z)^{3}}
$$
and therefore, we have that $g'(r_{\mathcal{S}}(\alpha))=0$. This observation shows that
$g$ cannot be univalent in $|z|<r$ if $r> r_{\mathcal{S}}(\alpha)$. Thus, $r_0=r_{\mathcal{S}}(\alpha)$.
\hfill $\Box$

\vspace{8pt}

\noindent \textbf{Proof of Theorem \ref{6-10th2}.}
Assume that $f,g\in \mathcal{U}$ with $\frac{f(z)+g(z)}{z}\neq 0$ in the unit disk $\ID$.
Then
$$\left|U_f(z)-1\right|\leq |z|^{2}
~\mbox{ and }~ \left|U_g(z)-1\right| \leq |z|^{2}.
$$
Further, let
$$\frac{z}{f(z)}=1+b_{1}z+b_{2}z^{2}+ \cdots ~\mbox { and }~
\frac{z}{g(z)}=1+c_{1}z+c_{2}z^{2}+\cdots ,
$$
where $b_1=-f''(0)/2$ and $c_1=-g''(0)/2$ so that $|b_1|\leq 2$ and $|c_1|\leq 2$.

Then, (\ref{8veq1a}) gives that
$$ \frac{z}{f(z)}=1 +b_1z-\int_0^1\frac{w_1(tz)}{t^2}\,dt ~\mbox{ and }~
 \frac{z}{g(z)}=1 +c_1z-\int_0^1\frac{w_2(tz)}{t^2}\,dt
,  \quad z\in \ID
$$
for some analytic functions $w_j\colon \ID \ra \ID$ such that $w_j(0)=w_j'(0)=0$ for $j=1,2$.
Thus,  for $|z|=r$, we have
\be\label{6-10eq5}
S_1=\frac{1}{2}\left|\frac{z}{f(z)}-\frac{z}{g(z)}\right|
\leq \frac{1}{2}|b_1-c_1|r+r^{2}\leq 2r+r^{2}
\ee
since $|b_1\pm c_1|\leq |b_1|+|c_1|\leq 4$. Also,  for $|z|=r$, we have that
$$
S_2=\frac{1}{2}\left|\frac{z}{f(z)}+\frac{z}{g(z)}\right|
\geq 1-\frac{|b_{1}+c_{1}|}{2}r-\sum_{n=2}^{\infty}\frac{|b_{n}|+|c_{n}|}{2}r^{n} .
$$
By the Cauchy-Schwarz inequality and Lemma \ref{5-10lem1}\textbf{(c)}, it follows that
$$\sum_{n=2}^{\infty}|b_{n}|r^{n} \leq \left(\sum_{n=2}^{\infty}(n-1)^{2}|b_{n}|^{2}\right)^{1/2}
\left(\sum_{n=2}^{\infty}\frac{r^{2n}}{(n-1)^{2}}\right)^{1/2}
\leq r\sqrt{A(r)}
$$
where
$$A(r)=\sum_{n=1}^{\infty}\frac{r^{2n}}{n^2}
\leq \left(\sum_{n=1}^{\infty}\frac{1}{n^4}\right)^{1/2}\left(\sum_{n=1}^{\infty}r^{4n}\right)^{1/2}
=\frac{\pi ^2}{\sqrt{90}} \frac{r^2}{\sqrt{1-r^4}}.
$$
Therefore, we have
$$\sum_{n=2}^{\infty}|b_{n}|r^{n}\leq r\sqrt{A(r)}\leq
\frac{\pi r^2}{\sqrt[4]{90} \,\sqrt[4]{1-r^4}}
$$
and similarly,
$$\sum_{n=2}^{\infty}|c_{n}|r^{n} \leq r\sqrt{A(r)}\leq
\frac{\pi r^2}{\sqrt[4]{90} \,\sqrt[4]{1-r^4}}.
$$
Using these two inequalities, we deduce that
\be\label{6-10eq6}
S_2\geq  1-2r- \frac{\pi r^{2} }{\sqrt[4]{90} \, \sqrt[4]{1-r^{4}}}.
\ee
Next we consider the function $F$ defined by (\ref{6-10eq4}) so that
$$U_F(z)-1 = \left(\frac{z}{F(z)}\right)^{2}F'(z)-1
=\left(\frac{2z}{f(z)+g(z)}\right)^{2}\frac{f'(z)+g'(z)}{2}-1.
$$
We may now rewrite the right hand expression and obtain
\beqq
|U_F(z)-1|&=&\left|2U_f(z) \frac{f^2(z)}{(f(z)+g(z))^{2}}
+2U_g(z)\frac{g^2(z)}{(f(z)+g(z))^{2}}
+\frac{(f(z)-g(z))^{2}}{(f(z)+g(z))^{2}}\right|\\
&\leq & \frac{2|z|^{2}(|f(z)|^{2}+|g(z)|^{2}) +|f(z)-g(z)|^{2}}{|f(z)+g(z)|^{2}}\\
&=&  \frac{|z|^{2}(|f(z)-g(z)|^{2}+|f(z)+g(z)|^{2}) +|f(z)-g(z)|^{2}}{|f(z)+g(z)|^{2}}\\
& = &  |z|^{2}+(1+|z|^{2})\frac{|f(z)-g(z)|^{2}}{|f(z)+g(z)|^{2}}\\
&=& |z|^{2}+(1+|z|^{2})\frac{S_1^2}{S_2^2}.
\eeqq
Using the estimates for $S_1$ and $S_2$ from (\ref{6-10eq5}) and (\ref{6-10eq6}), we find that for $|z|=r$,
$$|U_F(z)-1| \leq
r^{2}+(1+r^{2})\frac{(2r+r^{2})^{2}}{\left(1-2r-\frac{\pi}{\sqrt[4]{90}}\frac{r^{2}}{\sqrt[4]{1-r^{4}}}\right)^{2}}.
$$
Thus, $|U_F(z)-1|<1$ for $|z|=r<r_0$ (and hence, by the maximum modulus theorem, $F\in {\mathcal U}$ for $|z|<r_0$)
if $0<r<r_{0}$, where $ r_{0}\approx 0.262453 $ is the root of the equation (\ref{6-10eq8}).
\hfill $\Box$

\vspace{8pt}

\bcor\label{6-10cor3}
Let $f,g\in \mathcal{U}$ with $\frac{f(z)+g(z)}{z}\neq 0$ in the unit disk $\ID$,
$f''(0)+g''(0)=0$, and $F$ be defined by {\rm (\ref{6-10eq4})}.
Then $F\in {\mathcal U}$ for $|z|<r_0$, where $r_{0} \approx 0.3512$ is the root of the equation
$$
r^{2}+(1+r^{2})\frac{(2r+r^{2})^{2}}{\left(1-\frac{\pi}{\sqrt[4]{90}}\frac{r^{2}}{\sqrt[4]{1-r^{4}}}\right)^{2}}=1.
$$
\ecor
\bpf
The hypothesis that $f''(0)+g''(0)=0$ gives that $b_1+c_1=0$ and therefore, it suffices to
observe that the estimate for $S_2$ in the
proof of Theorem \ref{6-10th2} takes the form
$$S_2\geq  1- \frac{\pi r^{2} }{\sqrt[4]{90} \, \sqrt[4]{1-r^{4}}}
$$
and the rest of the proof is similar.
\epf

Similarly, by the obvious observation in the proof of Theorem \ref{6-10th2}, we easily have the following.

\bcor\label{6-10cor4}
Let $f,g\in \mathcal{U}$ with $\frac{f(z)+g(z)}{z}\neq 0$ in the unit disk $\ID$,
either $f''(0)=0$ or $g''(0)=0$, and $F$ be defined by {\rm (\ref{6-10eq4})}.
Then $F\in {\mathcal U}$ for $|z|<r_0$, where $r_{0}\approx 0.400502$ is the root of the equation
$$r^{2}+(1+r^{2})\frac{(r+r^{2})^{2}}{\left(1-r-\frac{\pi}{\sqrt[4]{90}}\frac{r^{2}}{\sqrt[4]{1-r^{4}}}\right)^{2}}=1.
$$
\ecor

\bcor\label{6-10cor5}
Let $f,g\in \mathcal{U}$ with $\frac{f(z)+g(z)}{z}\neq 0$ in the unit disk $\ID$,
$f''(0)=g''(0)=0$, and $F$ be defined by {\rm (\ref{6-10eq4})}.
Then $F\in {\mathcal U}$ for $|z|<r_0$, where $r_{0}\approx 0.667827$ is the root of the equation
$$r^{2}+(1+r^{2})\frac{r^4}{\left(1-\frac{\pi}{\sqrt[4]{90}}\frac{r^{2}}{\sqrt[4]{1-r^{4}}}\right)^{2}}=1.
$$
\ecor

As in Corollaries \ref{6-10cor3}-\ref{6-10cor5}, Theorem \ref{6-10th3} may be stated with an improved
form in the cases where either $f''(0)+g''(0)=0$, or $f''(0)=0$ or $g''(0)=0$, or $f''(0)=g''(0)=0$, respectively.
The sharpness of radii $r_0$ in Theorem \ref{6-10th2} and related corollaries are open.
Moreover, if
$$\mathcal{U}_2:=\{f\in \mathcal{U}:\, f''(0)=0\},
$$
then the inequality (\ref{6-10eq9}) shows that each $f\in \mathcal{U}_2$
satisfies the condition ${\rm Re\,}(f(z)/z)>1/2$ in $\ID$
and hence, $\mathcal{U}_2\subset \mathcal{N}$. Therefore, it is natural to ask
the analog of Corollary \ref{6-10cor5} when $f,g\in {\mathcal N}.$

\end{document}